\documentclass[a4paper,12pt]{article}
\usepackage[utf8]{inputenc}

\usepackage{amsmath}
\usepackage{amssymb,latexsym}
\usepackage{amsthm}

\usepackage[top=3cm, left=3cm, right=2cm, bottom=3cm]{geometry}
\newtheorem{theorem}{Theorem}
\newtheorem{lemma}[theorem]{Lemma}
\newtheorem{proposition}[theorem]{Proposition}

\title{Explicit bounds for solutions to  optimal control problems}
\author{Miguel Oliveira, Georgi Smirnov\\~\\
University of Minho}
\date{~}
\begin{document}

\maketitle

\begin{abstract}
In this paper we present explicit bounds for optimal control in a Lagrange problem without end-point constraints.
The approach we use is due to Gamkrelidze and is based on the equivalence of the Lagrange problem and a time-optimal problem for differential inclusions.

\end{abstract}

\section{Introduction}

The first works concerning regularity of solutions to basic problem in calculus of variations appeared more than a century ago  \cite{B,Tonelli}.
In the last 30 years, regularity of solution to problems of calculus of variations and optimal control has been a subject of intensive studies  
(see, e.g.,  \cite{ClV,ClV2,ClV3,Cellina,CF,CFM,Z}). However, at least to our knowledge, there are no explicit bounds for optimal control in Lagrange problems.
 In this paper, we obtain such  bounds. This is done under rather strong conditions. The approach we use is very close to the one from \cite{ST,T} and is based on the equivalence between the Lagrange problem  and a time-optimal control problem. The time-optimal control problem is reduced to a time-optimal problem for a differential inclusion. Application of 
 necessary conditions of optimality to this problem allows us to deduce explicit bounds for the control in the original Lagrange problem.

We consider the following Lagrange problem
\begin{eqnarray}
&& \int_0^1 L(t,x(t),u(t))dt \rightarrow\inf,\nonumber\\
&& \dot{x}(t)=g(t,x(t))u(t),\label{1}\\
&& x(0)=0, \nonumber
\end{eqnarray}
where $L:R^n\times R^m\rightarrow R$ and $g$  is an $n\times m$-matrix, and obtain explicit bounds for the optimal control in two special cases:
\begin{enumerate}
\item The autonomous case: $L=L(x,u)$ and $g=g(x)$;
\item The case with $g=g(t)$.
\end{enumerate}

  Problem (\ref{1}) is equivalent to the time optimal control problem
\begin{eqnarray}
&& T\rightarrow\inf,\nonumber\\
&& \frac{d(t,y)(\tau )}{d\tau}
=
\frac{(1,g(t(\tau), y(\tau))w(\tau ))}{L(t(\tau), y(\tau ),w(\tau ))},\label{2}\\
&& (t,y)(0)=(0,0),\;\; t(T)=1.\nonumber
\end{eqnarray}
The equivalence of (\ref{1}) and (\ref{2}) has been discussed by Gamkrelidze \cite{G} in the case of the basic problem of calculus of variations and by Sarychev and Torres in the optimal control framework \cite{ST}.

We denote the norm of the vector $x\in R^n$ by $|x|$. The closed unit ball in $R^n$ is denoted by $B_n$.  The distance between $x\in R^n$ and $C\subset R^n$ is denoted by $d(x,C)$. The tangent cone to $C$ at $x\in C$ is defined as $T(C,x)=\{ v\in R^n\mid \lim_{\lambda\downarrow 0}\lambda^{-1}d(x+\lambda v,C)=0\}$. The conjugate cone to a cone $K\subset R^n$ is denoted by $K^*=\{ x^*\in R^n \mid \langle x^*,x\rangle\geq 0,\; x\in K\}$. The set of absolutely continuous functions $x:[t_0,t_1]\rightarrow R^n$ is denoted by $AC([t_0,t_1],R^n)$.

\section{Main results}

Assume that $L(\cdot, \cdot,\cdot)$ and $g(\cdot, \cdot)$ are continuously differentiable functions satisfying the following conditions:

\vspace{5mm}

\noindent (C1) There exists a function $\theta:r\rightarrow R$ satisfying $L(t,x,u)\geq \theta (|u|)>0$, for all $t,x,u$, and $\lim_{r\rightarrow\infty}{r}/{\theta(r)}=0$.

\vspace{5mm}

\noindent (C2) Function $L(t,x,\cdot)$ is strictly convex, i.e., there exists a constant $\mu>0$ such that
$$
L(t,x,u)+\langle\nabla_u L(t,x,u),v-u\rangle +\frac{\mu}{2}|v-u|^2\leq L(t,x,v)
\;\;{\rm for\: all}\;\; t,x,u,v.
$$
\noindent (C3) The following growth condition is satisfied:
$$
|\nabla_{(t,x)}L(t,x,u)||g(t,x)u|\leq \xi L(t,x,u)+\delta \;\;{\rm for\: all}\;\; t,x,u,
$$
where $\xi>0$ and $\delta>0$ are some constants.

\noindent (C4) 
There exist constants $c_g>0$ and $c_{\nabla g}>0$ such that
$$
|g(t,x)|\leq c_g,\;\;\; |\nabla_{(t,x)} {g}(t,x)|\leq c_{\nabla g}\;\;{\rm for\: all}\;\; t,x.
$$ 

\vspace{5mm}

Obviously we have $\theta(r)/r\geq 1$, whenever $r\geq r_0$, if $r_0>0$ is big enough.
Set
$$
c=r_0+\int_0^1L(t,0,0)dt.
$$
Let $(\hat{u}(\cdot),\hat{x}(\cdot))\in L_1([0,1],R^m)\times AC([0,1],R^n)$ be the solution to Lagrange problem (\ref{1}).
 Denote by $M$ the set of points $t\in[0,1]$ where the inequality $|{\hat{u}}(t)|\leq \theta(|{\hat{u}}(t)|)$ holds.
Since
$$
 c\geq r_0+ \int_0^1 L(t,\hat{x}(t),{\hat{u}}(t))dt\geq r_0+ \int_M L(t,\hat{x}(t),{\hat{u}}(t))dt
$$
\begin{equation}
 \geq \int_{[0,1]\setminus M} r_0dt+\int_M\theta(|{\hat{u}}(t)|)dt\geq \int_0^1|{\hat{u}}(t)|dt,
\label{ner1}
\end{equation}
we have
$$
|\hat{x}(t)|\leq c_g c,
$$ 
whenever $t\in [0,1]$.
Set $\Omega=[0,1]\times c_g c B_n$, 
$$
\Lambda_0=\max_{(t,x)\in\Omega}L(t,x,0),
$$ 
$$
\Lambda_1=\max_{(t,x)\in\Omega}|\nabla_u L(t,x,0)|,
$$ 
and 
$$
\sigma (r) = \max_{(t,x)\in\Omega,\; |u|\leq r} (\langle\nabla_uL(t,x,u),u\rangle-L(t,x,u)).
$$ 
The following Lemma is an immediate consequence of (C2).
\begin{lemma}
\label{lem00}
The following inequality holds:
$$
-\Lambda_0-\Lambda_1|u|+\frac{\mu}{2}|u|^2\leq L(t,x,u).
$$
The function $\sigma (r)$ tends to infinity as $r\rightarrow\infty$.
\end{lemma}

{\em Proof}. The inequality is a consequence of the following:
$$
L(t,x,0)+\langle\nabla_u L(t,x,0),u\rangle + \frac{\mu}{2}|u|^2\leq L(t,x,u).
$$
Next, from condition (C2) we have
$$
L(t,x,u)-\langle\nabla_u L(t,x,u),u\rangle+\frac{\mu}{2}|u|^2\leq L(t,x,0).
$$
Therefore
$$
\frac{\mu}{2}r^2-\Lambda_0\leq \sigma (r).\;\;\; \qed
$$

\vspace{5mm}

From this lemma we see that there exists $0<T_0<1$ such that $\beta=\sigma((c+1)/T_0)>{\frac{\delta}{\xi}}$.

\subsection*{Autonomous case}

Assume that $L$ and $g$ do not depend on $t$.

\begin{theorem}
\label{th1}
The following inequality holds
$$
|{\hat{u}}(t)|\leq \ell=\max\left\{ \sqrt{\frac{2}{\mu}\left(\Lambda_0+\beta\right)},\frac{\Lambda_1+\sqrt{\Lambda_1^2+4\mu\Lambda_0}}{2}\right\}.
$$
\end{theorem}

\subsection*{Nonautonomous case}

Now assume that $g$ does not depend on $x$. Set 
$$
\eta=\sup_{r\geq 0}\frac{r}{\theta(r)+\beta}
$$
and
$$
\gamma=\frac{c_{\nabla g}+c_g\xi}{c_g\xi}e^{c_g\eta\xi (\Lambda_0+\beta)}.
$$

\begin{theorem}
\label{th2}
The following inequality holds
$$
|\hat{u}(t)|\leq  \ell=\max\left\{\sqrt{\frac{2}{\mu}(\Lambda_0+\beta)(1+\gamma \xi)},
\frac{\Lambda_1+\sqrt{\Lambda_1^2+4\mu\Lambda_0}}{2}
\right\}
$$
\end{theorem}

\section{Background notes}

We shall use the following proposition (cf. \cite{G,ST}).

\begin{proposition}
\label{pr1}
The following assertions hold:
\begin{enumerate}
\item For any admissible control process $(u(\cdot),x(\cdot ))$ of problem (\ref{1}) there exists a trajectory $(t,y)(\tau )$, $\tau\in [0,T]$ of control system (\ref{2}) such that
$(t,y)(0)=(t',x(t'))$, $(t,y)(T)=(t'',x(t''))$, and
\begin{equation}
\label{T1}
T=\int_{t'}^{t''} L(t,x(t),{u}(t))dt.
\end{equation}
\item For any trajectory $(t,y)(\tau )$, $\tau\in [0,T]$ of control system (\ref{2}) such that $\frac{d}{d\tau}(t,y)\neq (0,0)$ almost everywhere, there exists 
a control process $(u(\cdot),x(\cdot ))$ of problem (\ref{1}) such that $x(t(0))=y(0)$, $x(t(T))=y(T)$, and
\begin{equation}
\label{T2}
T=\int^{t(T)}_{t(0)} L(t,x(t),{u}(t))dt.
\end{equation}
\end{enumerate}
\end{proposition}

{\em Proof}. Let $(u(\cdot),x(\cdot ))$ be an admissible control process of problem (\ref{1}) and $t',t''\in [0,1]$, $t'<t''$. The function
$$
\tau (t)=\int_{t'}^t L(s,x(s),u(s))ds,\;\; t\in [t',1]
$$
is bounded, strictly monotonous and absolutely continuous. Its inverse, $t=t(\tau)$ is also strictly monotonous and absolutely continuous. Therefore the function
$y(\tau)=x(t(\tau))$ is absolutely continuous and the function $w(\tau)=u(t(\tau))$ is measurable. Moreover the equalities
$$
\frac{dt(\tau)}{d\tau}=\frac{1}{L(t(\tau),y(\tau),w(\tau))},
$$
$$
\frac{dy(\tau)}{d\tau}=\frac{g(t(\tau),y(\tau))w(\tau)}{L(t(\tau),y(\tau),w(\tau))},
$$
and (\ref{T1}) hold.

Now, let us consider a trajectory $(t,y)(\tau )$, $\tau\in [0,T]$ of control system (\ref{2}) such that $\frac{d}{d\tau}(t,y)\neq (0,0)$ almost everywhere. 
Since the function $t=t(\tau)$ is strictly monotonous and absolutely continuous, the inverse function $\tau=\tau (t)$ is  strictly monotonous and absolutely continuous. Hence the function $x(t)=y(\tau (t))$ is absolutely continuous and the function $u(t)=w(\tau(t))$ is measurable. Differentiating $x(t)$, we obtain
$$
\dot{x}(t)=\frac{dy}{d\tau}\frac{d\tau}{dt}=g(t,x(t))w(t).
$$
Since $d\tau/dt=L$, we get (\ref{T2}). $\qed$

\vspace{5mm}

Let $F:R^n\rightarrow R^n$ be a Lipschitzian set-valued map with compact values. Consider the following time-optimal problem
\begin{eqnarray*}
&& T\rightarrow\min,\\
&& \dot{x}\in F(x),\\
&& x(0)=x_0,\;\; x(T)\in S.
\end{eqnarray*}
Here $S\subset R^n$ is a convex set. 
Let $\hat{x}(\cdot)\in AC([0,T],R^n)$ be a solution to this problem. Consider a convex cone $K(t)\subset T({\rm gr\:co}F,(\hat{x}(t),\dot{\hat{x}}(t)))$ measurably depending on $t\in [0,T]$.

There exist many necessary conditions of optimality for time-optimal problems with differential inclusions (see e.g. \cite{Cl,V}).
For our considerations the most suitable formulation is contained in the  following proposition which is a consequence of  \cite[Theorem 5]{S1}.
\begin{proposition}
\label{pr2}
There exists a function $p(\cdot )\in AC([0,T],R^n)$ such that
\begin{enumerate}
\item $(\dot{p},p)\in -K^*(t)$, $\langle p(t),\dot{\hat{x}}(t)\rangle\equiv h\geq 0$;
\item $p(T)\in (T(S,\hat{x}(T)))^*$;
\item $|p(T)|>0$.
\end{enumerate}
\end{proposition}

In the case of a smooth control system, the Yorke approximation can be chosen as the cone $K(t)$.
Let $U\subset R^k$, and let $f:R^n\times U\rightarrow R^n$ be a function. Assume that $f$ is differentiable in $x$ and the set $f(x,U)$ is convex for all $x\in R^n$. For $(\hat{x},\hat{u})\in R^n\times U$ denote $\hat{v}=f(\hat{x},\hat{u})$ and set $C=\nabla_x f(\hat{x},\hat{u})$, $K=T(f(\hat{x},U),\hat{v})$. Recall the following proposition \cite[p. 38]{S2}.
\begin{proposition}
\label{pr3}
The following inclusion holds:
$$
\{ (x,v)\in R^n\times R^n\mid v\in Cx+K\}\subset T((\hat{x},\hat{u}), {\rm gr} f(\cdot,U)).
$$
\end{proposition}

Recall also the following useful formula \cite[p. 50]{S2}.
\begin{proposition}
\label{pr4}
Let $C:R^n\rightarrow R^n$ be a linear operator, and let $K\subset R^n$ be a convex cone. Then the following equality holds:
$$
\{ (x,v)\in R^n\times R^n\mid v\in Cx+K\}^*
$$
$$
=
\{ (x^*,v^*)\in R^n\times R^n\mid x^*=-C^*v^*,\; v^*\in K^*\}.
$$
\end{proposition}

\section{Auxiliary propositions}

Let  us consider the set-valued map
$$
G(t,y)=\left\{ (v^0,v)\in R\times R^n\mid 
v^0=\frac{\rho }{ L(t,y,q)+\beta},\right.
$$
$$
\left. v=\frac{\rho g(t,y) w}{ L(t,y,w)+\beta},\; w\in U,\; \rho\in [0,1]\right\}.
$$
\begin{lemma}\label{lem1}
The set-valued map $G$ has convex compact values and is Lipschitzian in $(t,y)$ in the set $\Omega$.
\end{lemma}

{\em Proof}. Let  $(v^0_i,v_i)\in G(t,y)$, $i=1,2$. There exist $\rho_i\in [0,1]$ and $w_i\in U$, $i=1,2$, such that
$$
v_i^0=\frac{\rho_i}{ L(t,y,w_i)+\beta},\;\;\;
 v_i=\frac{\rho_i g(t,y)w_i}{ L(t,y,w_i)+\beta},\;\;\; i=1,2.
$$
Let  $\alpha_1,\; \alpha_2 \geq 0$, $\alpha_1+\alpha_2 =1$. Show that 
 $\alpha_1 (v^0_1,v_1) +\alpha_2 (v^0_2,v_2)\in G(t,y)$. Put 
$$
\alpha'_1= \frac{\alpha_1\rho_1}{ L(t,y,w_1)+\beta}
\left( \frac{\alpha_1\rho_1}{ L(t,y,w_1)+\beta}
+\frac{\alpha_2\rho_2}{ L(t,y,w_2)+\beta}
\right)^{-1} 
$$
 and
$$
\alpha_2'= \frac{\alpha_2\rho_2}{ L(t,y,w_2)+\beta}
\left( \frac{\alpha_1\rho_1}{ L(t,y,w_1)+\beta}
+\frac{\alpha_2\rho_2}{L (t,y,w_2)+\beta}
\right)^{-1}. 
$$
Obviously $\alpha'_1,\; \alpha_2' \geq 0$, $\alpha'_1+\alpha_2' =1$. 
Set $w=\alpha_1'w_1+\alpha_2'w_2$. Then we get $ L(t,y,w)\leq\alpha_1'L(t,y,w_1)+\alpha_2' L(t,y,w_2)$.
We have
$$
\alpha_1 v_1+\alpha_2 v_2=
\frac{\alpha_1\rho_1 g(t,y)w_1}{ L(t,y,w_1)+\beta}+
\frac{\alpha_2\rho_2 g(t,y)w_2}{ L(t,y,w_2)+\beta}
$$
$$
=
\left(\frac{\alpha_1\rho_1}{ L(t,y,w_1)+\beta}+
\frac{\alpha_2\rho_2}{ L(t,y,w_2)+\beta}\right)g(t,y) w.
$$
Put
$$
\rho=
\left(\frac{\alpha_1\rho_1}{ L(t,y,w_1)+\beta}+
\frac{\alpha_2\rho_2}{ L(t,y,w_2)+\beta}\right)( L(t,y,w)+\beta).
$$
Then we obtain
$$
\rho\leq 
\frac{\alpha_1\rho_1}{ L(t,y,w_1)+\beta}( L(t,y,w_1)+\beta)+
\frac{\alpha_2\rho_2}{ L(t,y,w_2)+\beta}(  L(t,y,w_2)+\beta)
$$
$$
=\alpha_1\rho_1+\alpha_2\rho_2\leq 1,
$$
i.e., $\rho\in [0,1]$. Therefore $G(t,y)$ is convex.

From (C1) we have
\begin{equation}
\label{int_***}
|v^0|\leq
\frac{1}{ \theta (|w|)+\beta},\;\;\;
|v|\leq
\frac{c_g|w|}{ \theta (|w|)+\beta},\;\;\;
{\rm for\:all}\;\;\; (v^0,v)\in G(t,y).
\end{equation}
Let
$$
(v_k^0,v_k)=
\left(\frac{\rho_k }{ L(t,y,w_k)+\beta},\frac{\rho_k g(t,y)w_k}{ L(t,y,w_k)+\beta}\right),
$$
where $w_k\in R^n$ and $\rho_k\in [0,1]$, be a sequence converging to a point  
$(v^0_0,v_0)$. If the sequence $w_k$ is bounded, then, without loss of generality, the sequence $(w_k,\rho_k)$ converges. Passing to the limit we get
 $(v^0_0,v_0)\in G(t,y)$. If the sequence $w_k$ is unbounded, then there exists a subsequence converging to infinity. Without loss of generality  $w_k$ goes to infinity. From inequalities   (\ref{int_***}) we obtain $(w_k^0,w_k)\rightarrow (0,0)$. Hence  $(w_0^0,w_0)=(0,0)\in G(t,y)$. Thus   $G(t,y)$ is a closed set. From (\ref{int_***}) we see that it is bounded.

Let $(t_1,y_1)$ and $(t_2,y_2)$ be two points in $\Omega$. 
Let 
$$
 (v_1^0,v_1)=\frac{(\rho,\rho g(t_1,y_1)w)}{ L(t_1,y_1,w)+\beta}\in G(t_1,y_1).
$$
Consider the point
$$
 (v_2^0,v_2)=\frac{(\rho,\rho g(t_2,y_2)w)}{ L(t_2,y_2,w)+\beta}\in G(t_2,y_2).
$$
Since $\beta>\delta/\xi$, from (C3) we have
$$
|v^0_1-v_2^0|\leq\max_{\lambda\in [0,1],w}\left|\nabla_{(t,x)}\left(\frac{1}{( L(\lambda t_1+(1-\lambda )t_2,\lambda y_1+(1-\lambda )y_2,w)+\beta)}\right)\right| 
$$
$$
\times (|t_1-t_2|+|y_1-y_2|)
\leq\frac{\xi}{\beta} (|t_1-t_2|+|y_1-y_2|)
$$
and
$$
|v_1-v_2|\leq\max_{\lambda\in [0,1],w}\left|\nabla_{(t,x)}\left(\frac{g(\lambda t_1+(1-\lambda )t_2,\lambda y_1+(1-\lambda )y_2)w}{( L(\lambda t_1+(1-\lambda )t_2,\lambda y_1+(1-\lambda )y_2,w)+\beta)}\right)\right| 
$$
$$
\times (|t_1-t_2|+|y_1-y_2|)
\leq\left( \frac{c_{\nabla g}}{\beta}+\eta\xi c_g\right) (|t_1-t_2|+|y_1-y_2|) ,
$$
i.e. $G$  is Lipschitzian in $y$ in the set $\Omega$. $\qed$

\vspace{5mm}

Let $(\hat{u}(\cdot),\hat{x}(\cdot))$ be a solution to problem (\ref{1}). By the first part of Proposition \ref{pr1} there exists
a trajectory $(\hat{t},\hat{y})(\tau )$, $\tau\in [0,\hat{T}]$ of control system 
\begin{equation}\label{3}
\frac{d(t,y)(\tau )}{d\tau}
=
\frac{(1,g(t(\tau),y(\tau))w(\tau ))}{ L(t(\tau ),y(\tau ),w(\tau ))+\beta},\;\; w(\tau)\in R^n.
\end{equation}
 such that
$(\hat{t},\hat{y})(0)=(0,0)$, $\hat{t}(\hat{T})=1$, and
$$
\hat{T}=\int_{0}^{1} ( L(t,\hat{x}(t),{\hat{u}}(t))+\beta)dt.
$$
The control corresponding to $\hat{y}(\cdot)$ is denoted by $\hat{w}(\cdot)$. 

\begin{lemma}
\label{lem10}
There exists a nonzero function  $(q,p)(\cdot)\in AC([0,\hat{T}],R\times R^n)$ such that
\begin{eqnarray}
&& \frac{dq(\tau )}{d\tau}=-\frac{\langle {\hat{g}_t}(\tau)\hat{w}(\tau),p(\tau)\rangle}{L(\hat{t}(\tau ),\hat{y}(\tau ),\hat{w}(\tau ))+\beta}\nonumber\\
&&+\frac{ (q(\tau )+\langle {g}(\hat{t}(\tau),\hat{y}(\tau))\hat{w}(\tau),p(\tau )\rangle )L_t(\hat{t}(\tau ),\hat{y}(\tau ),\hat{w}(\tau )) }{( L(\hat{t}(\tau ),\hat{y}(\tau ),\hat{w}(\tau ))+\beta)^2},\label{P1}\\
&& \frac{dp(\tau )}{d\tau}= -\frac{\nabla_x\langle {g}(\hat{t}(\tau),\hat{y}(\tau))\hat{w}(\tau),p(\tau)\rangle}{L(\hat{t}(\tau ),\hat{y}(\tau ),\hat{w}(\tau ))+\beta}\nonumber\\
&&+\frac{ (q(\tau )+\langle {g}(\hat{t}(\tau),\hat{y}(\tau))\hat{w},p(\tau )\rangle )\nabla_x L(\hat{t}(\tau ),\hat{y}(\tau ),\hat{w}(\tau )) }{( L(\hat{t}(\tau ),\hat{y}(\tau ),\hat{w}(\tau ))+\beta)^2},\label{P2}\\
&& p(\hat{T})=0,\label{P2*}\\
&&\label{P3}
\frac{g^*(\hat{t}(\tau),\hat{y}(\tau))p(\tau )}{ L(\hat{t}(\tau ),\hat{y}(\tau ),\hat{w}(\tau ))+\beta}\nonumber\\
&&-\frac{(q(\tau )+\langle {g}(\hat{t}(\tau),\hat{y}(\tau))\hat{w}(\tau),p(\tau )\rangle )\nabla_w L(\hat{t}(\tau ),\hat{y}(\tau ),\hat{w}(\tau ))}{( L(\hat{t}(\tau ),\hat{y}(\tau ),\hat{w}(\tau ))+\beta)^2}=0,\\
&&\label{P4} \frac{ q(\tau )+\langle {g}(\hat{t}(\tau),\hat{y}(\tau))\hat{w}(\tau),p(\tau )\rangle  }{ L(\hat{t}(\tau ),\hat{y}(\tau ),\hat{w}(\tau ))+\beta}\equiv h>0.
\end{eqnarray}
\end{lemma}

{\em Proof}. From the second part of Proposition \ref{pr1} we see that  $(\hat{t},\hat{y})(\tau )$, $\tau\in [0,\hat{T}]$ is a solution to the problem 
\begin{eqnarray}
&& T\rightarrow\inf,\nonumber\\
&& \frac{d(t,y)(\tau )}{d\tau}
=
\frac{(1,{g}(t(\tau),y(\tau))w(\tau ))}{ L(t(\tau ),y(\tau ),w(\tau ))+\beta},\;\; w(\tau)\in R^n, \label{4}\\
&& (t,y)(0)=(0,0),\;\;\; t(T )=1.\nonumber
\end{eqnarray}
The time-optimal problem
\begin{eqnarray}
&& T\rightarrow\inf,\nonumber\\
&& \frac{d(t,y)(\tau )}{d\tau}\in G(t,y),\label{TG}\\
&& (t,y)(0)=(0,0),\;\;\; t(T )=1\nonumber,
\end{eqnarray}
also has a solution $(\tilde{t},\tilde{y})(\tau)$, $\tau\in [0,\tilde{T}]$. By the Filippov lemma there exists a measurable function $(\tilde{\rho},\tilde{w}) (\tau)$, $\tau\in [0,\tilde{T}]$, such that
$$
\frac{d(\tilde{t},\tilde{y})(\tau )}{d\tau}
=
\frac{\tilde{\rho}(\tau)(1,{g}(t(\tau))\tilde{w}(\tau ))}{ L(\tilde{t}(\tau ),\tilde{y}(\tau ),\tilde{w}(\tau ))+\beta}
$$
at almost all points where $d(\tilde{t},\tilde{y})/d\tau\neq (0,0)$. 
Applying Propositions \ref{pr2}-\ref{pr4}, we see that there exist  $(q,p)(\cdot)\in AC([0,\tilde{\tau}],R\times R^n)$, a nonzero function, and a constant $h\geq 0$ such that
\begin{eqnarray}
&& \frac{dq(\tau )}{d\tau}=-\frac{\tilde{\rho}(\tau)\langle {\hat{g}_t}(\tau)\tilde{w}(\tau),p(\tau)\rangle}{L(\tilde{t}(\tau ),\tilde{y}(\tau ),\tilde{w}(\tau ))+\beta}\nonumber\\
&&+\frac{ \tilde{\rho}(\tau)(q(\tau )+\langle {g}(\tilde{t}(\tau),\tilde{y}(\tau))\tilde{w}(\tau),p(\tau )\rangle )L_t(\tilde{t}(\tau ),\tilde{y}(\tau ),\tilde{w}(\tau )) }{( L(\tilde{t}(\tau ),\tilde{y}(\tau ),\tilde{w}(\tau ))+\beta)^2},\label{n1}\\
&& \frac{dp(\tau )}{d\tau}= -\frac{\tilde{\rho}(\tau)\nabla_x\langle {g}(\tilde{t}(\tau),\tilde{y}(\tau))\tilde{w}(\tau),p(\tau)\rangle}{L(\tilde{t}(\tau ),\tilde{y}(\tau ),\tilde{w}(\tau ))+\beta}\nonumber\\
&&+\frac{ \tilde{\rho}(\tau)(q(\tau )+\langle {g}(\tilde{t}(\tau),\tilde{y}(\tau))\tilde{w},p(\tau )\rangle )\nabla_x L(\tilde{t}(\tau ),\tilde{y}(\tau ),\tilde{w}(\tau )) }{( L(\tilde{t}(\tau ),\tilde{y}(\tau ),\tilde{w}(\tau ))+\beta)^2},\label{n2}\\
&& p(\tilde{T})=0,\label{n2*}\\
&& h\equiv\frac{\tilde{\rho}(\tau)(q(\tau)+\langle {g}(\tilde{t}(\tau),\tilde{y}(\tau))\tilde{w}(\tau),p(\tau)\rangle )}{ L(\tilde{t}(\tau ),\tilde{y}(\tau ),\tilde{w}(\tau ))+\beta}\nonumber\\
&&\geq
\frac{{\rho}(q(\tau)+\langle {g}(\tilde{t}(\tau),\tilde{y}(\tau)){w},p(\tau)\rangle )}{ L(\tilde{t}(\tau ),\tilde{y}(\tau ),{w})+\beta},\;\;\rho\in [0,1],\;\; w\in R^n,\label{n3}
\end{eqnarray}
at almost all points such that $d(\tilde{t},\tilde{y})/d\tau\neq (0,0)$. 

Let us show that $d(\tilde{t},\tilde{y})/d\tau\neq (0,0)$ almost everywhere. 
If $d(\tilde{t},\tilde{y})/d\tau = (0,0)$ on a set of positive measure, then $h=0$. At points where $\tilde{\rho} (\tau)>0$, from maximum condition (\ref{n3}) we have
$$
\frac{\tilde{\rho}(\tau)g^*(\tilde{t}(\tau),\tilde{y}(\tau))p(\tau )}{ L(\tilde{t}(\tau ),\tilde{y}(\tau ),\tilde{w}(\tau ))+\beta}
$$
\begin{equation}
\label{n31}
-\frac{\tilde{\rho}(\tau)(q(\tau )+\langle {g}(\tilde{t}(\tau),\tilde{y}(\tau))\tilde{w}(\tau),p(\tau )\rangle )\nabla_w L(\tilde{t}(\tau ),\tilde{y}(\tau ),\tilde{w}(\tau ))}{( L(\tilde{t}(\tau ),\tilde{y}(\tau ),\tilde{w}(\tau ))+\beta)^2}=0
\end{equation}
Since $q(\tau)+\langle {g}(\tilde{t}(\tau),\tilde{y}(\tau))\tilde{w}(\tau),p(\tau)\rangle =0$, from (\ref{n2}) and (\ref{n2*}) we get $p=0$ and, hence, $q=0$, a contradiction. Thus $\tilde{\rho}(\tau)=0$ almost everywhere. This is impossible. Hence $d(\tilde{t},\tilde{y})/d\tau\neq (0,0)$ at almost all points $\tau\in [0,\tilde{T}]$. 

Therefore conditions (\ref{n1})-(\ref{n3}) are satisfied almost everywhere and $\tilde{\rho}(\tau)>0$ at almost all points $\tau\in [0,\tilde{T}]$. Thus $h>0$, because the equality $h=0$ implies, as above, that $(q,p)(\tau)\equiv 0$. From   (\ref{n3}) we obtain $\tilde{\rho}(\tau)=1$. Thus we can identify the trajectories $(\hat{t},\hat{y})(\cdot)$ and $(\tilde{t},\tilde{y})(\cdot)$. Both of them are solutions to time-optimal problem (\ref{TG}) and satisfy necessary conditions of optimality (\ref{n1}), (\ref{n2}), and (\ref{n31}) with $\tilde{\rho}=1$. $\qed$

\vspace{5mm}
Denote by $\hat{\tau}(\cdot)$ the function inverse to $\hat{t}(\cdot)$. Then we have 
${\hat{u}}(\cdot)=\hat{w}(\hat{\tau}(\cdot))$. Therefore it suffices to obtain the bounds for $\hat{w}(\cdot)$.
 We shall use the notations $\hat{L}(\tau)$ for $L(\hat{t}(\tau),\hat{y}(\tau),\hat{w}(\tau))$ and $\hat{g}(\tau)$ for $g(\hat{t}(\tau),\hat{y}(\tau))$.

\begin{lemma}
\label{lem2}
If $q(\tau)\leq 0$, then $|\hat{w}(\tau)|>( c+1)/T_0$.
\end{lemma}

{\em Proof}. Multiplying (\ref{P3}) by $ \hat{w}(\tau)$, we obtain
$$
( \hat{L}(\tau)+\beta)\langle \hat{g}(\tau)\hat{w}(\tau),p(\tau)\rangle
$$
\begin{equation}
\label{P44}
= (q(\tau)+\langle  \hat{g}(\tau)\hat{w}(\tau),p(\tau)\rangle)\langle\nabla_w\hat{L}(\tau), \hat{w}(\tau)\rangle.
\end{equation}
Since $q(\tau)\leq 0$, we have $\langle  g(\hat{t}(\tau))\hat{w}(\tau),p(\tau)\rangle>0$. From (\ref{P44}) we get
$$
 \langle\nabla_w\hat{L}(\tau),\hat{w}(\tau)\rangle=\frac{\langle  g(\hat{t}(\tau))\hat{w}(\tau),p(\tau)\rangle}{q(\tau)+\langle  g(\hat{t}(\tau))\hat{w}(\tau),p(\tau)\rangle}( \hat{L}(\tau)+\beta)\geq  \hat{L}(\tau)+\beta.
$$
From this we obtain
$$
\beta\leq \langle\nabla_u\hat{L}(\tau),\hat{w}(\tau)\rangle-\hat{L}(\tau)\leq \sigma (|\hat{w}(\tau)|).
$$
Hence,  we have
$$
|\hat{w}(\tau)|\geq\sigma^{-1}(\beta)= \sigma^{-1}\left(\sigma\left(\frac{ c+1}{T_0}\right)\right)=\frac{c+1}{T_0}.\;\;\;\qed
$$

\section{Proof of Theorem \ref{th1}}

If $L$ and $g$ do not depend on $t$, then from (\ref{P1}) we have $dq/d\tau=0$. Combining Lemma \ref{lem2} with (\ref{ner1}), we see that $q$ is a positive constant. 
From condition (C2) we have
$$
\hat{L}(\tau)-\langle \nabla_w\hat{L}(\tau),\hat{w}(\tau)\rangle +\frac{\mu}{2}|\hat{w}(\tau)|^2\leq \Lambda_0.
$$
From this and (\ref{P44}) we obtain
$$
\frac{\mu}{2}|\hat{w}(\tau)|^2\leq \Lambda_0-\hat{L}(\tau)+\langle \nabla_w\hat{L}(\tau),\hat{w}(\tau)\rangle=
$$
$$
\Lambda_0-\hat{L}(\tau)+\frac{\langle p(\tau),g(\hat{y}(\tau))\hat{w}(\tau)\rangle}{(q(\tau)+\langle p(\tau),g(\hat{y}(\tau))\hat{w}(\tau)\rangle)}(\hat{L}(\tau)+\beta).
$$
If $\langle p(\tau),g(\hat{y}(\tau))\hat{w}(\tau)\rangle> 0$, then we have
$$
\frac{\mu}{2}|\hat{w}(\tau)|^2\leq \Lambda_0-\hat{L}(\tau)+\hat{L}(\tau)+\beta=\Lambda_0+\beta.
$$
Hence
\begin{equation}
\label{nn1}
|\hat{w}(\tau)|\leq\sqrt{\frac{2}{\mu}\left( \Lambda_0+\beta\right)}.
\end{equation}
If $\langle p(\tau),\hat{w}(\tau)\rangle\leq 0$, then from Lemma \ref{lem00} we get
$$
\frac{\mu}{2}|\hat{w}(\tau)|^2\leq \Lambda_0-\hat{L}(\tau)\leq \Lambda_0+\Lambda_1 |\hat{w}(\tau)| -\frac{\mu}{2} |\hat{w}(\tau)|^2.
$$
Thus we obtain 
\begin{equation}
\label{nn2}
|\hat{w}(\tau)|\leq \frac{\Lambda_1+\sqrt{\Lambda_1^2+4\mu\Lambda_0}}{2}.\;\;\;\qed
\end{equation}

\section{Proof of Theorem \ref{th2}}

If $q(\tau)\geq 0$, then as in the proof of Theorem \ref{th1} we get
$$
|{\hat{w}}(\tau)|\leq =\max\left\{ \sqrt{\frac{2}{\mu}\left(\Lambda_0+\beta\right)},\frac{\Lambda_1+\sqrt{\Lambda_1^2+4\mu\Lambda_0}}{2}\right\}.
$$
Show that if $q(\tau_1)= 0$ and $q(\tau)< 0$, $\tau\in ]\tau_1,\tau_2]$, then $|q(\tau)|/|p(\tau)|\leq {\gamma}$, $\tau\in [\tau_1,\tau_2]$.
Since $q(\tau)< 0$, $\tau\in ]\tau_1,\tau_2]$, we have
$$
\frac{q(\tau)+\langle p(\tau),\hat{w}(\tau)\rangle}{|p(\tau)|( \hat{L}(\tau)+\beta)}\leq \frac{|\hat{w}(\tau)|}{\theta(|\hat{w}(\tau)|)+\beta}\leq\eta.
$$
From (\ref{P1}), (\ref{P2}), and condition (C3) we get
\begin{eqnarray}
&& \left|\frac{dq(\tau)}{d\tau}\right| \leq \frac{|\hat{g}_t(\tau)||\hat{w}(\tau)||p(\tau)|}{\hat{L}(\tau)+\beta} +\frac{ (q(\tau)+\langle p(\tau),\hat{g}(\tau)\hat{w}(\tau)\rangle)|\hat{L}_t(\tau)|}{( \hat{L}(\tau)+\beta)^2}\nonumber\\
&&\leq  c_{\nabla g}\eta |p(\tau)|+c_g \eta |p(\tau)| \frac{ |\hat{L}_t(\tau)|}{ \hat{L}(\tau)+\beta}\leq
(c_{\nabla g}\eta +c_g \eta\xi) |p(\tau)| ,\label{u1}\\
&& \left|\frac{dp(\tau)}{d\tau}\right| \leq \frac{ (q(\tau)+\langle p(\tau),\hat{g}(\tau)\hat{w}(\tau)\rangle)|\nabla_x\hat{L}(\tau)|}{( \hat{L}(\tau)+\beta)^2}\nonumber\\
&&\leq c_g\eta |p(\tau)| \frac{ |\nabla_x\hat{L}(\tau)|}{ \hat{L}(\tau)+\beta}\leq c_g\eta\xi|p(\tau)|.\label{u2}
\end{eqnarray}
whenever $\tau\in [\tau_1,\tau_2]$. 
From this we obtain
$$
\frac{d}{d\tau}\frac{|q(\tau)|}{|p(\tau)|}\leq\frac{|dq(\tau)/d\tau||p(\tau)|+|q(\tau)||dp(\tau)/d\tau|}{|p(\tau)|^2}
$$
$$
\leq\eta\left( c_{\nabla g}+c_g\xi+c_g\xi\frac{|q(\tau)|}{|p(\tau)|}\right).
$$
Since $q(\tau_1)=0$, applying the Gronwall inequality we have
\begin{equation}
\label{utt}
\frac{|q(\tau)|}{|p(\tau)|}\leq \frac{c_{\nabla g}+c_g\xi}{c_g\xi} e^{c_g\eta\xi (\tau_2-\tau_1)}.
\end{equation}
Observe that $(\bar{t},\bar{y})(\cdot)$ where $\bar{t}(\cdot)$ is the solution to the differential equation
$$
\frac{dt}{d\tau}=\frac{1}{L(t,0,0)+\beta}
$$
satisfying $\bar{t}(0)=0$, and $\bar{y}\equiv 0$, is an admissible solution to the time-optimal problem
(\ref{TG}) on an interval $[0,\bar{T}]$, where $\bar{T}$ is such that $\bar{t}(\bar{T})=1$. From the inequality
$$
\frac{1}{L(t,0,0)+\beta}\geq \frac{1}{\Lambda_0+\beta}
$$
we get $\hat{T}\leq \bar{T}\leq \Lambda_0+\beta$. Therefore we have $\tau_2-\tau_1\leq \hat{T}\leq \Lambda_0+\beta$.
Thus, from (\ref{utt}) we obtain $|q(\tau)|\leq\gamma |p(\tau)|$.

Let us rewrite (\ref{P44}) in the form
$$
\langle\hat{g}(\tau)\hat{w}(\tau),p(\tau)\rangle=h\langle\hat{L}_w(\tau),\hat{w}(\tau)\rangle
$$
Combining this with (\ref{P4}) we get
$$
h=\frac{q(\tau)+h\langle\hat{L}_w(\tau),\hat{w}(\tau)\rangle}{\hat{L}(\tau)+\beta}
$$
Hence from (C2) we have
$$
\frac{q(\tau)}{h}=\hat{L}(\tau)+\beta- \langle\hat{L}_w(\tau),\hat{w}(\tau)\rangle\leq\Lambda_0-\frac{\mu}{2}|\hat{w}(\tau)|^2+\beta.
$$
Thus for $\tau\in [\tau_1,\tau_2]$ we get
$$
\frac{\mu}{2}|\hat{w}(\tau)|^2\leq \Lambda_0+\beta+\frac{|q(\tau)|}{h}\leq \Lambda_0+\beta+\frac{\gamma |p(\tau)|}{h}.
$$
To evaluate $|p(\tau)|/h$, observe that (\ref{P2}) and (C3) imply
$$
\frac{d|p(\tau)|/h}{d\tau}\leq \frac{|\nabla_x \hat{L}(\tau)|}{\hat{L}(\tau)+\beta}\leq\xi.
$$
Since $p(\hat{T})=0$, we obtain $|p(\tau)|/h\leq \hat{T}\xi\leq (\Lambda_0+\beta)\xi$. This ends the proof. $\qed$

\section*{Acknowledgements}
The authors are grateful to Delfim Torres for bibliographical support. This work was partially supported by project PTDC/EEI-AUT/2933/2014 (TOCCATA), funded by Project 3599 - Promover a 
Produ\c c\~ao Cient\'\i fica e Desenvolvimento Tecnol\'ogico
e a Constitui\c c\~ao de Redes Tem\'aticas (3599-PPCDT) and FEDER funds through
COMPETE 2020, Programa Operacional Competitividade e Internacionaliza\c c\~ao (POCI),
and by national  funds through Funda\c c\~ao para a Ci\^encia e a Tecnologia (FCT). The work of Miguel Oliveira was supported by FCT through the PhD fellowship SFRH/BD/111854/2015.

\end{document}